\documentclass[ECP,preprint,NODS]{ejpecp}
\usepackage{xcolor}
\usepackage[T1]{fontenc}
\usepackage[utf8]{inputenc}
\usepackage{amsmath,amssymb,amsthm}
\usepackage[authoryear]{natbib}
\usepackage{chngcntr}

\SHORTTITLE{Removing the floor from Hoeffding and Bernstein inequalities}

\TITLE{Removing the floor from Hoeffding and Bernstein\\ inequalities for finite-population U-statistics}

\AUTHORS{
  Jianhang~Ai\footnote{Faculty of Electrical Engineering, Czech Technical University in Prague, Prague, Czech Republic.   \EMAIL{aijianha@fel.cvut.cz}}
\and 
  Yiming~Chen\footnote{School of Mathematical Sciences, Peking University, China.   \EMAIL{ymchenmath@math.pku.edu.cn}}
\and 
Ondřej Kuželka\footnote{Faculty of Electrical Engineering, Czech Technical University in Prague, Prague, Czech Republic.\EMAIL{ondrej.kuzelka@fel.cvut.cz}}
  \and 
  Yuyi Wang\footnote{CRRC Zhuzhou Institute \& Tengen Intelligence Institute,
Zhuzhou, China. \EMAIL{yuyiwang920@gmail.com}}
}


\AMSSUBJ{60E15, 62D05, 60C05}

\SUBMITTED{}

\ABSTRACT{
Let $C$ be a uniformly sampled $n$-element subset of a finite population $D$ of
size $N$, and let $U_C$ be the complete U-statistic of order $k$ built from a
bounded kernel on the $k$-subsets of $D$. We prove a product inequality
for the uniform distribution on $\binom{D}{n}$ and deduce from it the comparison
$\log\E e^{\theta(U_C-\E U_C)}\le (n/k)\log\E e^{(k\theta/n)(Y-\E Y)}$, where $Y$
is a single kernel value taken at a uniformly sampled $k$-subset of $D$. This
yields Hoeffding, Bennett and Bernstein inequalities for $U_C$ with effective
sample size $n/k$ in place of the $\lfloor n/k\rfloor$ available so far. The gain
is negligible when $n$ is much larger than $k$, but it approaches a factor of two
in the exponent when $k$ is close to $n$.
}

\counterwithout{equation}{section}

\counterwithout{theorem}{section}

\theoremstyle{ejpecpbodyrm}
\newtheorem{example}{Example}

\newcommand{\E}{\mathbb{E}}
\newcommand{\Prob}{\mathbb{P}}
\newcommand{\Var}{\operatorname{Var}}
\newcommand{\R}{\mathbb{R}}
\newcommand{\Nat}{\mathbb{N}}

\begin{document}

\section{Introduction}
\label{sec:intro}

U-statistics were introduced by \citet{hoeffding1948}, who also proved the first
exponential concentration inequalities for them \citep{hoeffding1963}. The
classical argument starts from the following observation. A complete U-statistic
of order $k$ computed from a sample of size $n$ can be written as an average of
averages: we split the sample into \emph{disjoint} blocks of size $k$, we average
the $\lfloor n/k\rfloor$ kernel evaluations on these blocks, and then we average
over the ways of splitting. The evaluations inside one block-average are
independent, so Jensen's inequality reduces any convex functional of the
U-statistic to the same functional of an average of $\lfloor n/k\rfloor$
independent copies of a single kernel value. The effective sample size that we
get out of the reduction is therefore $\lfloor n/k\rfloor$ and not $n/k$.

For complete U-statistics of independent observations, \citet{wang2017} noticed
that this integer loss is unnecessary: instead of matching disjoint blocks
integrally, they used a fractional matching of the complete $k$-uniform
hypergraph on the sample, whose fractional matching number is exactly $n/k$. For
the finite-population setting, \citet{ai2022} proved Hoeffding- and
Bernstein-type inequalities with effective sample size $\lfloor n/k\rfloor$, and
\citet{ai2023} proved bounds that also reflect the size of the population.
\citet{blanc2025} proved a convex-function comparison between $U_C$ and an
average of $\lfloor n/k\rfloor$ independent kernel evaluations. This comparison
also recovers the inequalities of \citet{ai2022}.

The purpose of this note is to show that the floor in these bounds is not
intrinsic to finite-population U-statistics. Writing $Y$ for the value of the
kernel at a uniformly selected $k$-subset of the full population, we prove that
\begin{equation}
\label{eq:main-intro}
  \log \E \exp\{\theta(U_C-\E U_C)\}
  \;\le\;
  \frac{n}{k}\,
  \log \E \exp\Big\{\frac{k\theta}{n}\,(Y-\E Y)\Big\}.
\end{equation}
This holds for every $\theta\in\R$ and is the without-replacement analogue of
the floor-free refinement of \citet{wang2017}. We deduce it from a product
inequality for the uniform distribution on $\binom{D}{n}$
(Theorem~\ref{thm:product}), whose proof is short and elementary.

The product inequality retains the fractional value $n/k$ under sampling
without replacement. Thus the floor in the earlier bounds reflects the
disjoint-block argument rather than the underlying dependence. Replacing
$\lfloor n/k\rfloor$ by $n/k$ usually gives a modest numerical improvement,
but it yields a cleaner bound that may be easier to apply.

Equation~\eqref{eq:main-intro} controls exponential moments only; it does
not give a floor-free version of Blanc's convex-order reduction. We return to
this distinction in Section~\ref{sec:remarks}.

Let $h$ denote the kernel, and suppose that $a\le h(T)\le b$ for every
$k$-subset $T$. For $t>0$, the one-sided Hoeffding bounds are floor-based
\citep{ai2022} and floor-free (this paper):
\[
\begin{aligned}
\Prob(U_C-\E U_C\ge t)
&\le \exp\!\left\{-\frac{2\lfloor n/k\rfloor t^2}{(b-a)^2}\right\}
&&\text{(with the floor)},\\
\Prob(U_C-\E U_C\ge t)
&\le \exp\!\left\{-\frac{2(n/k)t^2}{(b-a)^2}\right\}
&&\text{(without the floor)}.
\end{aligned}
\]
The lower-tail formulas are identical; the two-sided bounds have an additional
factor of $2$.

\begin{example}
Consider a $[0,1]$-valued kernel with $k=10$, $n=169$, and $t=0.3$.
Substitution gives $0.056$ for the floor-based bound and $0.048$ for the
floor-free bound. Numerically this is modest, but it instantiates the structural
statement above: the bound depends on $n/k$, not on an integer block count.
\end{example}


\section{Setup}
\label{sec:setup}

Let $D$ be a finite set with $|D|=N$, and for an integer $r$ let
$\binom{D}{r}=\{S\subseteq D:|S|=r\}$. Throughout, we use the convention
$\binom{m}{j}=0$ for $j<0$ and for $j>m$. Fix integers $1\le k\le n\le N$, let
$a\le b$ and let $h:\binom{D}{k}\to[a,b]$ be a bounded kernel. Let $C$ be
uniformly distributed on $\binom{D}{n}$ and define the complete U-statistic
\[
  U_C=\frac{1}{\binom{n}{k}}\sum_{T\in\binom{C}{k}} h(T).
\]
Let $T_\star$ be uniformly distributed on $\binom{D}{k}$ and put $Y=h(T_\star)$,
$\mu=\E Y$ and $\sigma^2=\Var(Y)$. Thus $Y$ is a single kernel value evaluated on
a uniformly sampled $k$-subset of the whole population, and not of the sample.
Every fixed $T\in\binom{D}{k}$ is contained in $C$ with probability
$\binom{N-k}{n-k}/\binom{N}{n}=\binom{n}{k}/\binom{N}{k}$, and therefore
\begin{equation}
\label{eq:unbiased}
  \E U_C=\mu .
\end{equation}
Finally, let $\psi(s)=\log\E e^{s(Y-\mu)}$. Since $Y$ is bounded, $\psi$ is
finite and infinitely differentiable on $\R$, with $\psi(0)=\psi'(0)=0$ and
$\psi''(0)=\sigma^2$.

\section{A product inequality for uniform sampling of fixed-size subsets}
\label{sec:product}

\begin{theorem}
\label{thm:product}
Let $D$ be a finite set with $|D|=N$ and let $1\le k\le n\le N$. For each
$T\in\binom{D}{k}$ let $q_T\ge 0$, and put $\omega=\binom{n-1}{k-1}$. If $C$ is
uniformly distributed on $\binom{D}{n}$, then
\begin{equation}
\label{eq:product}
  \E_C \prod_{T\in\binom{C}{k}} q_T^{1/\omega}
  \;\le\;
  \Bigg(\frac{1}{\binom{N}{k}}\sum_{T\in\binom{D}{k}} q_T\Bigg)^{\!n/k},
\end{equation}
where the expectation is with respect to the uniformly sampled set $C$.
\end{theorem}

\begin{remark}
We believe that Theorem~\ref{thm:product} admits a sharper finite population version. For \(1\leq k\leq n<N\), we expect that \eqref{eq:product} can be strengthened to
\[
    \mathbb{E}_{C}
    \prod_{T\in\binom{C}{k}}
    q_T^{\,m_{N,n,k}/\binom{n}{k}}
    \leq
    \left(
        \frac{1}{\binom{N}{k}}
        \sum_{T\in\binom{D}{k}} q_T
    \right)^{m_{N,n,k}} .
\]
where $m_{N,n,k}=
    \frac{\log\!\bigl(N/(N-n)\bigr)}
         {\log\!\bigl(N/(N-k)\bigr)}.$
Equivalently, the effective sample size \(n/k\) in Corollary~\ref{cor:laplace} should be
replaceable by \(m_{N,n,k}\). Note that \(m_{N,n,k}>n/k\) whenever
\(k<n<N\). We do not pursue this refinement in the
present paper.
\end{remark}

Before we prove the theorem, we make two observations about it. First, both sides
of \eqref{eq:product} are homogeneous of degree $n/k$: the left-hand side is a
product of $\binom{n}{k}$ factors, each raised to the power $1/\omega$, and
$\binom{n}{k}/\omega=n/k$. The weights $1/\omega$ are exactly the optimal
fractional matching of the complete $k$-uniform hypergraph on $n$ vertices, which
is what \citet{wang2017} used in the independent case. The theorem shows
that these weights remain valid under sampling without replacement.

Second, the proof is an induction in which we delete a single element $x$ from
$D$. Conditionally on $x\notin C$, the induction hypothesis applies directly.
Conditionally on $x\in C$, the product splits into the factors that avoid $x$ and
the factors that use $x$, and the identity
$\binom{n-1}{k-1}=\binom{n-2}{k-1}+\binom{n-2}{k-2}$ says exactly that the two
groups of exponents form a pair of conjugate Hölder weights. The fractional
exponent then comes from Bernoulli's inequality, whose exponent is the real
number $n/k$ and need not be an integer.

\begin{proof}
We argue by induction on $N$, simultaneously for all pairs $(k,n)$ with
$1\le k\le n\le N$. If $N=1$, then $k=n=1$, and both sides of
\eqref{eq:product} equal $q_D$. Write
\[
  \bar q=\binom{N}{k}^{-1}\sum_{T\in\binom{D}{k}}q_T .
\]
Then \eqref{eq:product} becomes
$\E\prod_{T\in\binom{C}{k}}q_T^{1/\omega}\le \bar q^{\,n/k}$.

We first handle the two boundary cases $k=n$ and $n=N$. If
$k=n$, then $\omega=1$ and $\binom{C}{k}=\{C\}$, so both sides of
\eqref{eq:product} equal $\bar q$. If $N=n$, then $C=D$, and the
arithmetic--geometric mean inequality gives
\[
  \prod_{T\in\binom{D}{k}} q_T^{1/\omega}
  =\Bigg(\prod_{T\in\binom{D}{k}} q_T^{1/\binom{N}{k}}\Bigg)^{\binom{N}{k}/\omega}
  \le \bar q^{\,\binom{N}{k}/\omega}
  =\bar q^{\,N/k},
\]
where the last equality holds because $\omega=\binom{N-1}{k-1}$ in this case and
$\binom{N}{k}/\binom{N-1}{k-1}=N/k$.

It remains to consider $k<n<N$. Assume inductively that
\eqref{eq:product} holds for population size $N-1$ and every admissible pair of
parameters. Fix $x\in D$ and put
$D_0=D\setminus\{x\}$, so that $|D_0|=N-1\ge n$. Let
\[
  A=\frac{1}{\binom{N-1}{k}}\sum_{T\in\binom{D_0}{k}} q_T,
  \qquad
  B=\frac{1}{\binom{N-1}{k-1}}\sum_{R\in\binom{D_0}{k-1}} q_{R\cup\{x\}} .
\]
Thus $A$ averages $q_T$ over the $k$-subsets of $D$ that avoid $x$,
whereas $B$ averages $q_T$ over those that contain $x$. Since
$\binom{N-1}{k}/\binom{N}{k}=(N-k)/N$ and
$\binom{N-1}{k-1}/\binom{N}{k}=k/N$,
\begin{equation}
\label{eq:avg-split}
  \bar q=\frac{N-k}{N}\,A+\frac{k}{N}\,B ,
\end{equation}
We will also need $\Prob(x\in C)=n/N$.

Conditionally on $x\notin C$, the set $C$ is uniformly distributed on
$\binom{D_0}{n}$, so the induction hypothesis on $D_0$ with the parameters
$(k,n)$ gives
\begin{equation}
\label{eq:cond-out}
  \E\Bigg[\prod_{T\in\binom{C}{k}} q_T^{1/\omega}\;\Bigg|\;x\notin C\Bigg]
  \le A^{n/k}.
\end{equation}
Conditionally on $x\in C$, write $C=C_0\cup\{x\}$ with $C_0$ uniformly
distributed on $\binom{D_0}{n-1}$ and split the product as
\[
  \prod_{T\in\binom{C}{k}} q_T^{1/\omega}
  = F\cdot G,
  \qquad
  F=\prod_{T\in\binom{C_0}{k}} q_T^{1/\omega},
  \qquad
  G=\prod_{R\in\binom{C_0}{k-1}} q_{R\cup\{x\}}^{1/\omega}.
\]
Put $\omega_0=\binom{n-2}{k-1}$ and $\omega_1=\binom{n-2}{k-2}$. These are the
values of $\omega$ that correspond to the sample size $n-1$ and to the orders $k$
and $k-1$, and it is easy to check that $\omega_0+\omega_1=\omega$. We claim that
\begin{equation}
\label{eq:cond-in}
  \E[FG]\le A^{(n-k)/k}B .
\end{equation}
If $k=1$, then $\omega=\omega_0=1$ and $\omega_1=0$, the second product consists
of the single deterministic factor $q_{\{x\}}=B$, and the induction hypothesis on
$D_0$ with the parameters $(1,n-1)$ gives $\E[FG]=B\,\E F\le BA^{n-1}$, which is
\eqref{eq:cond-in}. If $k\ge 2$, then $\omega_0\ge 1$ and $\omega_1\ge 1$, so
Hölder's inequality with the conjugate exponents $\omega/\omega_0$ and
$\omega/\omega_1$ applies and yields
\[
  \E[FG]\le
  \big(\E F^{\omega/\omega_0}\big)^{\omega_0/\omega}
  \big(\E G^{\omega/\omega_1}\big)^{\omega_1/\omega}.
\]
Here $F^{\omega/\omega_0}=\prod_{T\in\binom{C_0}{k}}q_T^{1/\omega_0}$ and
$G^{\omega/\omega_1}=\prod_{R\in\binom{C_0}{k-1}}q_{R\cup\{x\}}^{1/\omega_1}$.
We can therefore apply the induction hypothesis on $D_0$ twice: first to
$(q_T)_{T\in\binom{D_0}{k}}$ with parameters $(k,n-1)$, and then to
$(q_{R\cup\{x\}})_{R\in\binom{D_0}{k-1}}$ with parameters $(k-1,n-1)$. This gives
$\E F^{\omega/\omega_0}\le A^{(n-1)/k}$ and
$\E G^{\omega/\omega_1}\le B^{(n-1)/(k-1)}$. Since
$\omega_0/\omega=(n-k)/(n-1)$ and $\omega_1/\omega=(k-1)/(n-1)$, we obtain
\[
  \E[FG]\le
  A^{\frac{n-1}{k}\cdot\frac{n-k}{n-1}}\,
  B^{\frac{n-1}{k-1}\cdot\frac{k-1}{n-1}}
  =A^{(n-k)/k}B,
\]
which is \eqref{eq:cond-in}. Combining \eqref{eq:cond-out} and
\eqref{eq:cond-in} with $\Prob(x\in C)=n/N$,
\begin{equation}
\label{eq:combined}
  \E\prod_{T\in\binom{C}{k}} q_T^{1/\omega}
  \le \frac{N-n}{N}\,A^{n/k}+\frac{n}{N}\,A^{(n-k)/k}B .
\end{equation}
If $A=0$, the right-hand side of \eqref{eq:combined} vanishes because $n>k$, and
there is nothing left to prove. If $A>0$, put $u=B/A\ge 0$ and divide
\eqref{eq:combined} by $A^{n/k}$; by \eqref{eq:avg-split} it suffices to show
that
\[
  \frac{N-n}{N}+\frac{n}{N}\,u
  \;\le\;
  \Big(\frac{N-k}{N}+\frac{k}{N}\,u\Big)^{\!n/k}.
\]
This is exactly Bernoulli's inequality $(1+z)^\alpha\ge 1+\alpha z$ with
$\alpha=n/k>1$ and $z=k(u-1)/N$. Indeed, the base on the right-hand side is
$1+z$, we have $\alpha z=n(u-1)/N$, so the left-hand side is $1+\alpha z$, and
$z\ge -k/N>-1$ because $u\ge0$ and $k<N$. This finishes the induction and the
proof.
\end{proof}

\begin{corollary}
\label{cor:laplace}
Recall that $\psi(s)=\log\E e^{s(Y-\mu)}$. Then, for every $\theta\in\R$,
\begin{equation}
\label{eq:laplace}
  \log \E \exp\{\theta(U_C-\mu)\}
  \;\le\;
  \frac{n}{k}\,\psi\Big(\frac{k\theta}{n}\Big).
\end{equation}
\end{corollary}

\begin{proof}
We apply Theorem~\ref{thm:product} with
$q_T=\exp\{\tfrac{k\theta}{n}(h(T)-\mu)\}$. Since
$\sum_{T\in\binom{C}{k}}(h(T)-\mu)=\binom{n}{k}(U_C-\mu)$ and
$\binom{n}{k}/\omega=n/k$, the left-hand side of \eqref{eq:product} becomes
$\E e^{\theta(U_C-\mu)}$ and the right-hand side becomes
$\big(\E e^{(k\theta/n)(Y-\mu)}\big)^{n/k}$. Taking logarithms then gives
\eqref{eq:laplace}.
\end{proof}

\begin{remark}
\label{rem:integer}
If $m=n/k$ is an integer and $Y_1,\dots,Y_m$ are independent copies of $Y$, the
right-hand side of \eqref{eq:laplace} is exactly the logarithmic moment
generating function of $m^{-1}\sum_{i=1}^m(Y_i-\mu)$.
Corollary~\ref{cor:laplace} extends this comparison to noninteger $n/k$.
\end{remark}

\begin{remark}
\label{rem:exact}
If $q_T=q>0$ for every $T$, both sides of \eqref{eq:product} equal
$q^{n/k}$. So \eqref{eq:product} is an equality in this case, and the exponent
$n/k$ cannot be changed. Indeed, replacing it by some $\alpha>n/k$ would fail
when $q<1$, whereas replacing it by some $\alpha<n/k$ would fail when $q>1$.
Note that this says
something about the exponent only, and nothing about the numerical constants in
Section~\ref{sec:conc}.
\end{remark}

\section{Concentration inequalities}
\label{sec:conc}

Corollary~\ref{cor:laplace} turns any bound on $\psi$ into a bound for $U_C$ with
$n/k$ in the role of a sample size. We now derive the standard consequences. We
put
\[
  c_+=b-\mu,\qquad c_-=\mu-a,\qquad c=\max\{c_+,c_-\},
\]
and we note that $c_++c_-=b-a$ and therefore $c\le b-a$.
If $\sigma^2=0$, then $h$ is constant on $\binom{D}{k}$, $U_C=\mu$ almost
surely, and all the tail probabilities below are zero. We therefore assume
$\sigma^2>0$ below. Then $Y$ is not constant, so $a<\mu<b$ and both $c_+$
and $c_-$ are strictly positive.

\begin{theorem}[Hoeffding inequality]
\label{thm:hoeffding}
For every $\theta\in\R$ we have
$\E \exp\{\theta(U_C-\mu)\}\le \exp\{k(b-a)^2\theta^2/(8n)\}$, and if $b>a$ then
\begin{equation}
\label{eq:hoeffding-tail}
  \Prob\big(|U_C-\mu|\ge t\big)
  \le 2\exp\Big\{-\frac{2n}{k}\cdot\frac{t^2}{(b-a)^2}\Big\},
  \qquad t>0 .
\end{equation}
\end{theorem}

\begin{proof}
Hoeffding's lemma, applied to $Y-\mu$, which takes values in an interval of
length $b-a$, gives $\psi(\lambda)\le \lambda^2(b-a)^2/8$ for all
$\lambda\in\R$. We substitute this into \eqref{eq:laplace} with
$\lambda=k\theta/n$ and get
$\tfrac{n}{k}\cdot\tfrac{k^2\theta^2}{n^2}\cdot\tfrac{(b-a)^2}{8}
=k(b-a)^2\theta^2/(8n)$. The tail bound \eqref{eq:hoeffding-tail} then follows
from the Chernoff bound applied to each of the two tails, optimized at
$\theta=\pm 4nt/(k(b-a)^2)$.
\end{proof}

\begin{corollary}
\label{cor:var}
The variance of the U-statistic satisfies
\[
  \Var(U_C)\le \frac{k}{n}\sigma^2 .
\]
\end{corollary}

\begin{proof}
Let $\varphi(\theta)=\log\E e^{\theta(U_C-\mu)}$ and
$F(\theta)=\tfrac{n}{k}\psi(k\theta/n)-\varphi(\theta)$. Both $\psi$ and
$\varphi$ are infinitely differentiable, because $Y$ and $U_C$ are bounded.
Corollary~\ref{cor:laplace} says that $F\ge 0$, and we have $F(0)=0$, so $0$ is a
global minimum of $F$ and hence $F''(0)\ge 0$. Since $Y-\mu$ and $U_C-\mu$ are
centered by \eqref{eq:unbiased}, we have $\psi''(0)=\sigma^2$ and
$\varphi''(0)=\Var(U_C)$, and therefore
$F''(0)=\tfrac{k}{n}\sigma^2-\Var(U_C)$, which is what we needed.
\end{proof}

In the remaining bounds we use the function $H(u)=(1+u)\log(1+u)-u$, $u\ge0$.

\begin{theorem}[Bennett and Bernstein inequalities]
\label{thm:bennett}
Assume $\sigma^2>0$ and let $t>0$. Then
\begin{equation}
\label{eq:bennett}
  \Prob(U_C-\mu\ge t)
  \le \exp\Big\{-\frac{n\sigma^2}{k\,c_+^2}\,
        H\Big(\frac{c_+ t}{\sigma^2}\Big)\Big\},
\end{equation}
and the same bound holds for $\Prob(\mu-U_C\ge t)$ with $c_-$ in place of $c_+$.
Consequently
\begin{equation}
\label{eq:bernstein}
  \Prob\big(|U_C-\mu|\ge t\big)
  \le 2\exp\Big\{-\frac{n}{k}\cdot\frac{t^2}{2(\sigma^2+ct/3)}\Big\}.
\end{equation}
\end{theorem}

\begin{proof}
Bennett's lemma \citep{bennett1962} states that a centered random variable $X$
with $X\le d$ almost surely, $d>0$, satisfies
\[
  \log\E e^{\lambda X}\le \frac{\Var(X)}{d^{2}}\big(e^{\lambda d}-1-\lambda d\big),
  \qquad \lambda\ge 0;
\]
see also the proof of Theorem 2.9 in \citet{blm2013}. We apply it to $X=Y-\mu$
with $d=c_+$ and substitute the result into \eqref{eq:laplace}. Writing
$\theta=(n/k)\lambda$ with $\lambda\ge0$, we get from the Chernoff bound that
\[
  \Prob(U_C-\mu\ge t)
  \le \inf_{\lambda\ge 0}
  \exp\Big\{\frac{n}{k}\Big(-\lambda t+\frac{\sigma^2}{c_+^2}
  \big(e^{\lambda c_+}-1-\lambda c_+\big)\Big)\Big\}.
\]
The exponent is convex in $\lambda$ and its derivative vanishes when
$e^{\lambda c_+}=1+c_+t/\sigma^2$. We substitute $u=c_+t/\sigma^2$, so that
$\lambda=c_+^{-1}\log(1+u)$ and $t=\sigma^2u/c_+$, and the exponent becomes
\[
  \frac{n}{k}\cdot\frac{\sigma^2}{c_+^2}
  \big(-u\log(1+u)+u-\log(1+u)\big)
  =-\frac{n\sigma^2}{k\,c_+^2}\,H(u),
\]
This is \eqref{eq:bennett}. If we apply the same argument to the kernel $-h$,
which interchanges $c_+$ and $c_-$ and leaves $\sigma^2$ unchanged, we get the
bound for the lower tail. Finally, \eqref{eq:bernstein} follows from
$H(u)\ge u^2/(2(1+u/3))$, valid for $u\ge0$, together with $c_\pm\le c$.
\end{proof}

\begin{corollary}
\label{cor:conf}
Assume $\sigma^2>0$ and let $\delta\in(0,1)$. With probability at least
$1-\delta$,
\begin{equation}
\label{eq:conf}
  |U_C-\mu|
  \le \sqrt{\frac{2k\sigma^2}{n}\log\frac{2}{\delta}}
      +\frac{kc}{3n}\log\frac{2}{\delta}.
\end{equation}
\end{corollary}

\begin{proof}
This does not follow from \eqref{eq:bernstein}, but from the moment bound behind
it. We combine Bennett's lemma with the elementary inequality
$e^x-1-x\le (x^2/2)/(1-x/3)$, valid for $0\le x<3$, and use $c_+\le c$; this
gives $\psi(\lambda)\le \lambda^2\sigma^2/(2(1-\lambda c/3))$ for
$0\le\lambda<3/c$. Substituting $\lambda=k\theta/n$ into \eqref{eq:laplace}, we
obtain
\[
  \begin{aligned}
    \log\E e^{\theta(U_C-\mu)}
      &\le \frac{\theta^2 v}{2(1-\theta c_1)},\\
    0\le\theta&<1/c_1,\qquad
    v=\frac{k\sigma^2}{n},\quad c_1=\frac{kc}{3n}.
  \end{aligned}
\]
In other words, $U_C-\mu$ is sub-gamma with variance factor $v$ and scale
parameter $c_1$. The same holds for $\mu-U_C$, by the argument with $-h$ used
above. The standard tail bound for sub-gamma random variables (see Section 2.4 of
\citealt{blm2013}) gives $\Prob(\pm(U_C-\mu)\ge\sqrt{2vs}+c_1s)\le e^{-s}$ for
every $s>0$. We take $s=\log(2/\delta)$ and add the two bounds, which finishes
the proof.
\end{proof}

\section{Remarks}
\label{sec:remarks}

On setting $c=b-a$, the bounds of Section~\ref{sec:conc} recover those of
\citet{ai2022} with $K=\lfloor n/k\rfloor$ replaced by $n/k$; retaining
$c=\max\{b-\mu,\mu-a\}$ gives a modest additional improvement in the range
term. Their Bernstein exponent, for instance, is
$-Kt^2/\big(2(\sigma^2+(b-a)t/3)\big)$. The floor had already been removed for
U-statistics of independent observations by \citet{wang2017}. What
Theorem~\ref{thm:product} adds is a way of handling the dependence that sampling
a subset of fixed cardinality creates.

Let us now compare Corollary~\ref{cor:laplace} with the reduction theorem of
\citet{blanc2025}. Let $p=n/k$, $q=\lfloor p\rfloor$, and let $Y_1,\dots,Y_q$ be
independent copies of $Y$. Blanc's theorem gives
$\E f(U_C)\le\E f(q^{-1}\sum_{i\le q}Y_i)$ for every convex $f$, and hence, if we
take $f(x)=e^{\theta x}$,
\begin{equation}
\label{eq:blanc}
  \log\E e^{\theta(U_C-\mu)}\le q\,\psi(\theta/q).
\end{equation}
Since $\psi$ is convex with $\psi(0)=0$, the slope function $s\mapsto\psi(s)/s$
is nondecreasing on $\R\setminus\{0\}$. Writing
$r\psi(\theta/r)=\theta\cdot\psi(\theta/r)/(\theta/r)$, we can see that
$r\mapsto r\psi(\theta/r)$ is nonincreasing on $(0,\infty)$ for every fixed
$\theta$: when $\theta<0$, the reversal caused by multiplication by $\theta$
compensates for the fact that $\theta/r$ increases with $r$. Since $p\ge q$,
Corollary~\ref{cor:laplace} is therefore at least as strong as
\eqref{eq:blanc}. It is strictly stronger for every $\theta\neq0$ when $Y$ is
nondegenerate and $p\notin\Nat$. The two statements are nevertheless of
different types, and ours is of the weaker type. Blanc's is a comparison in the
convex order with an actual random variable, and it controls all convex test
functions. Ours is about exponential moments only. For noninteger $p$ we do
not see any natural random variable that would play the role of an average of
$p$ independent copies of $Y$: the function
$s\mapsto(\E e^{s(Y-\mu)/p})^{p}$ need not be a moment generating function at
all. We also do not know whether the convex-order comparison itself can be made
floor-free.

Finally, we mention three limitations. First, our bounds carry no
finite-population correction and they are not sharp when $n/N$ is close to one.
For $n=N$ we have $\Var(U_C)=0$, whereas Corollary~\ref{cor:var} gives
$k\sigma^2/N$. The reason is that Theorem~\ref{thm:product} degenerates into the
arithmetic--geometric mean inequality when $n=N$, and a floor-free version of the
population-size-dependent bounds of \citet{ai2023} would therefore require a
different argument. Second, Remark~\ref{rem:exact} shows that the exponent $n/k$
is exact, but we make no claim about the constants $1/8$, $1/3$ and $2$ in
Section~\ref{sec:conc}. These are inherited from the corresponding inequalities
for a single kernel value. Third, \eqref{eq:conf} involves the population
variance $\sigma^2$, which is typically unknown. An empirical version, in the
spirit of \citet{maurer2009}, would be more directly usable.

\section*{Acknowledgement of AI use}

GPT-5.6 Sol and Claude Fable were used iteratively for brainstorming and
text editing, but the main ideas came from the human authors---in particular, the idea and techniques for removing the floor function follows YW's line of work, which previously lead to a similar case in the sampling with replacement setting \citep{wang2017}. All AI-assisted material was reviewed by the authors, who take
full responsibility for the manuscript.


\begin{thebibliography}{99}
\footnotesize

\bibitem[Ai et~al.(2022)]{ai2022}
Ai, J., Ku\v{z}elka, O., Wang, Y., 2022.
Hoeffding and Bernstein inequalities for U-statistics without replacement.
Statist. Probab. Lett. 187, 109528.

\bibitem[Ai et~al.(2023)]{ai2023}
Ai, J., Ku\v{z}elka, O., Wang, Y., 2023.
Hoeffding--Serfling inequality for U-statistics without replacement.
J. Theoret. Probab. 36, 390--408.

\bibitem[Bennett(1962)]{bennett1962}
Bennett, G., 1962.
Probability inequalities for the sum of independent random variables.
J. Amer. Statist. Assoc. 57, 33--45.

\bibitem[Blanc(2025)]{blanc2025}
Blanc, G., 2025.
Subsampling suffices for adaptive data analysis.
J. ACM 72 (1), Article 2.

\bibitem[Boucheron et~al.(2013)]{blm2013}
Boucheron, S., Lugosi, G., Massart, P., 2013.
Concentration Inequalities: A Nonasymptotic Theory of Independence.
Oxford University Press, Oxford.

\bibitem[Hoeffding(1948)]{hoeffding1948}
Hoeffding, W., 1948.
A class of statistics with asymptotically normal distribution.
Ann. Math. Statist. 19, 293--325.

\bibitem[Hoeffding(1963)]{hoeffding1963}
Hoeffding, W., 1963.
Probability inequalities for sums of bounded random variables.
J. Amer. Statist. Assoc. 58, 13--30.

\bibitem[Maurer and Pontil(2009)]{maurer2009}
Maurer, A., Pontil, M., 2009.
Empirical Bernstein bounds and sample variance penalization.
In: Proceedings of the 22nd Conference on Learning Theory (COLT).

\bibitem[Wang et~al.(2017)]{wang2017}
Wang, Y., Guo, Z.-C., Ramon, J., 2017.
Learning from networked examples.
Proc. Mach. Learn. Res. 76, 641--666.

\end{thebibliography}
\end{document}